\documentclass[sn-mathphys]{sn-jnl}

\usepackage{amssymb,amsmath,mathrsfs, amsthm}
\usepackage{hyperref}

\usepackage{cases}

\jyear{2021}%
\theoremstyle{thmstyleone}%
\numberwithin{equation}{section}
\newtheorem{theorem}{Theorem}[section]
\newtheorem{question}{Question}[section]%

\newtheorem{corollary}{Corollary}[section]
\theoremstyle{thmstyletwo}%
\newtheorem{remark}[theorem]{Remark}%
\newtheorem{lemma}{Lemma}[section]%
\newtheorem{conjecture}{Conjecture}[section]%

\theoremstyle{thmstylethree}%
\newtheorem{definition}[theorem]{Definition}%

\begin{document}

\title[Rigidity of the K\"ahler--Finsler Kobayashi metric and K\"ahler bases]{K\"ahler--Finsler rigidity of the Kobayashi metric and the K\"ahler base of weakly K\"ahler--Finsler metrics}



\author[]{\fnm{Xiaoshu} \sur{Ge}}
\affil{\orgdiv{Institute of Mathematics, Academy of Mathematics and Systems Science}, \orgname{Chinese Academy of Sciences}, \city{Beijing}, \postcode{100190}, \country{China}\\
E-mail: gxs@amss.ac.cn}




\abstract{
In this paper, we establish a rigidity theorem: for every balanced, smooth, bounded pseudoconvex domain, its Kobayashi metric is a smooth, strongly pseudoconvex, K\"ahler--Finsler metric restricted on a neighborhood of the origin if and only if the domain is a Hermitian ellipsoid, biholomorphic to the unit ball.
Thus for any balanced, smooth, bounded, strictly convex, and strongly pseudoconvex domain that is not biholomorphic to the unit ball, its Kobayashi metric is a non-K\"ahler--Finsler, weakly K\"ahler--Finsler metric. We also prove that the base manifold of every weak K\"ahler--Finsler metric is K\"ahler.
}




\maketitle

\section{Introduction and main results}

In complex Finsler geometry, there are four natural notions for generalizations of the K\"ahler condition (\cite{AP,GZ,G}):
\begin{itemize}
  \item strongly K\"ahler--Finsler: $\Gamma^i_{j;k}-\Gamma^i_{k;j}=0;$
  \item K\"ahler--Finsler: $(\Gamma^i_{j;k}-\Gamma^i_{k;j})v^j=0;$
  \item weakly K\"ahler--Finsler: $G_i(\Gamma^i_{j;k}-\Gamma^i_{k;j})v^j=0;$
\item new (weakly K\"ahler--Finsler) condition: $G_i(\Gamma^i_{j;k}-\Gamma^i_{k;j})=0.$
\end{itemize}
The first three are more well-known, while the last one proposed very recently by the author and X. Zhou in their manuscript \cite{GZ}.
In \cite{CS}, B.Chen and Y.Shen proved that,
the strongly K\"ahler--Finsler condition coincides with the K\"ahler--Finsler condition; and very recently, the author proved the coincidence for the last two in \cite{G}, and provide smooth non-K\"ahler--Finsler, weakly K\"ahler--Finsler examples\footnotemark[1]\footnotetext[1]{In \cite{CGZhou}, N. Cui, J. Guo and L. Zhou also provided  some non-K\"ahler--Finsler, weakly K\"ahler--Finsler examples more early, but their examples have singularity with respect the fiber variables besides the zero section. From a K\"ahler-Berwald rigidity theorem of X. Ge and X. Zhou in \cite{GZ}, together with the coincidence results for K\"ahler--Finsler and K\"ahler--Berwald in complex dimension two by N. Aldea and G. Munteanu in \cite{AM}, the Kobayashi metric of smooth bounded strongly convex domains of complex dimension two provide such kind examples, which indeed contain the special examples in \cite{G}. However, the examples in \cite{G} were checked when the authors of \cite{G,GZ} do not realize the results of \cite{AM}. Zisu Zhao also checked some special examples of complex dimension two recently.}.
Thus the study about the difference between strongly K\"ahler--Finsler and weakly K\"ahler--Finsler conditions is an important question.

It is known that every K\"ahler--Finsler metric has a K\"ahler base manifold (cf. \cite{Aikou,CSZ}). Notice that all known weakly K\"ahler--Finsler metrics have K\"ahler base, one may wonder whether there exists weakly K\"ahler--Finsler examples on non-K\"ahler manifolds. Such examples will become non-K\"ahler--Finsler weakly K\"ahler--Finsler examples automatically,
and make weakly K\"ahler--Finsler manifolds a larger class, in manifold sense, then K\"ahler manifolds as a generalization concept.
This induces a nature question:
\begin{question}\label{ques}
Wether the base manifold for any K\"ahler--Finsler metric is K\"ahler?
\end{question}
In this paper, using the stronger form for the weakly K\"ahler--Finsler condition, we obtain the positive answer to Question \ref{ques}.
\begin{theorem}[cf. Theorem \ref{pf}]\label{thm}
 Every K\"ahler--Finsler metric has a K\"ahler base.
\end{theorem}

Notice that the only known smooth weakly K\"ahler--Finsler metrics are provided by the Kobayashi metrics of smooth bounded strongly convex domains(see \cite{Lem81,Lem82,Pang,AP2,Pat}), which have many good properties, such as the constance holomorphic sectional curvature $-4$. In complex dimension two, they are always non-K\"ahler--Finsler whenever the domains are not biholomorphic to the unit ball $\mathbb B^2$ \cite{G,GZ}. For higher dimension, whether they are non-K\"ahler--Finsler is still a question, and they are conjectured to be whenever the domains are not biholomorphic to $\mathbb B^n$ in \cite{GZ}. This Kobayashi rigidity question is another main theme of the present paper.

There is another motivation for considering the Kobayashi rigidity theorem, induced from a well known conjecture of  S.-Y. Cheng in \cite{Cheng79}, which conjectured that:
the Bergman metric of a smooth bounded strongly pseudoconvex domain $D\in\mathbb C^n$ is K\"ahler-Einstein if and only if $D$ is biholomorphic to $\mathbb B^n$ (more generally, Yau conjectured the homogeneity for the pseudoconvex case \cite{Yau}). In 2021, X. Huang and M. Xiao conformed Cheng's conjecture (\cite{Huang-Xiao}).
Consequently, H. Gaussier and A. Zimmer proved that: the Kobayashi metric of a $C^2$-smooth bounded strongly pseudoconvex domain $D\subset \mathbb C^n$ is K\"ahler if and only if the universal cover of $D$ is biholomorphic to $\mathbb B^n$ (\cite{Gau-Zim}); R. Dong, R. Wang and B. Wong obtained the parallel results for the Carath\'eodory metric: the Carath\'eodory metric of a $C^2$-smooth bounded strongly pseudoconvex domain $D\subset \mathbb C^n$ is K\"ahler if and only if $D$ is biholomorphic to $\mathbb B^n$ (\cite{DWW}).

The author and X. Zhou weakened the K\"ahler condition in \cite{Gau-Zim,DWW} to the K\"ahler--Berwald condition in \cite{GZ},
since K\"ahler--Finsler implies K\"ahler--Berwald in complex dimension two by results of N. Aldea and G. Munteanu in \cite{AM},  they established the following rigidity result.
\begin{corollary}[\cite{GZ}]\label{KF-rigidity}
Let $D\subset\mathbb C^2$ be a ($C^2$-)smooth, bounded, strongly pseudoconvex domain. Then
$F_K$ is K\"ahler--Finsler if and only if the universal cover of $D$ is biholomorphic to $\mathbb B^2$.
If $D$ is in addition strongly convex, then $F_K$ is K\"ahler--Finsler if and only if $D$ is biholomorphic to $\mathbb B^2$.
\end{corollary}
Thus they obtained: for every ($C^2$-)smooth, bounded, strongly convex domain $D\subset\mathbb C^2$, $F_K$ is non-K\"ahler--Finsler, weakly K\"ahler--Finsler, unless $D$ is biholomorphic to the unit ball. Then they proposed the following conjecture.
\begin{conjecture}[a special version, \cite{GZ}]\label{conj}
For every smooth bounded strongly convex (resp. pseudoconvex) domain $D\subset \mathbb C^n,n\geq2$,
the Kobayashi metric is K\"ahler--Finsler if and only if $D$ (resp. the universal cover of $D$) is biholomorphic to the unit ball $\mathbb B^n.$
\end{conjecture}

Based on the author's check of the Kobayashi metrics for a special class of smooth, bounded, strongly convex domains of complex dimension two in \cite{G}, we obtain a rigidity theorem for balanced, smooth, bounded, strongly convex domains of general dimensions.
\begin{theorem}[cf. Theorem \ref{rigi-KF}]\label{Kob-KF-rig}
For every balanced, smooth, bounded, strongly convex domain
$D\subset\mathbb C^n$, its Kobayashi metric is a smooth strongly pseudoconvex K\"ahler--Finsler metric
if and only if $D$ is a Hermitian ellipsoid, linearly biholomorphic to the unit ball $\mathbb B^n$.
\end{theorem}
\begin{remark}
We can replace the balanced condition by circular, since they coincide under the strongly convex setting.
The proof of Theorem \ref{Kob-KF-rig} is independent of the results obtained in \cite{DWW,Gau-Zim,GZ,Huang-Xiao}.
\end{remark}

Besides convexity condition for the domains, M.-Y. Pang obtained that any balanced, smooth, bounded, strongly pseudoconvex domain $D$ admits a neighborhood $U\subset D$ of the origin, such that for any  $x\in U$, there is a Monge--Amp\`ere exhaustion $\tau^x$ with center $x$;
G. Patrizio and A. Spiro's obtained the smooth dependence of the Monge--Amp\`ere exhaustions away from the center of smooth, bounded, strongly pseudoconvex domains.
Together their results provide the smoothness of the Kobayashi metric in a neighborhood of the origin for any balanced, smooth, bounded, strongly pseudoconvex domain, which generalizing the strictly convex Lempert theory. Then we can generalize Theorem \ref{Kob-KF-rig} to the strongly pseudoconvex setting.
And for balanced pseudoconvex domain, the smooth, strongly pseudoconvex, K\"ahler--Finsler conditions of $F_K$ near the origin can help us weaken the boundary regularity condition to pseudoconvex.
Hence we establish the following rigidity theprem, which conforms conjecture \ref{conj} for the balanced case.
\begin{theorem}[cf. Theorem \ref{thm:balanced-rigidity}]\label{rig-main}
Let $D\subset\mathbb C^n$, $n\geq2$, be a balanced, bounded, pseudoconvex domain. Then its Kobayashi metric is a smooth, strongly pseudoconvex, K\"ahler--Finsler metric restricted on a neighborhood of the origin if and only if $D$ is a Hermitian ellipsoid, linearly biholomorphic to $\mathbb B^n$.
\end{theorem}
\begin{remark}
Theorem \ref{rig-main} conforms Conjecture \ref{conj} under the balanced condition.
\end{remark}

Notice that the Kobayashi metric $F_K$ of every smooth, strongly convex domain is smooth (see \cite{Lem81,Lem82,Pang}), and the existence of holomorphic geodesic disks makes $F_K$ weakly K\"ahler--Finsler when $F_K$ is strongly pseudoconvex (see \cite{AP2})\footnotemark[1]\footnotetext[1]{In \cite{AP2}, the weakly K\"ahler--Finsler condition was called K\"ahler there.}.
Thus a direct corollary of Theorem \ref{rig-main} provides the following results.
\begin{corollary}
Let $D\subset \mathbb C^n$ be any balanced, smooth, bounded, strongly convex domain. Then the Kobayashi metric of $D$ is non-K\"ahler--Finsler, weakly K\"ahler--Finsler unless $D$ is biholomorphic to the unit ball.
\end{corollary}

The rest of this paper is organized as follows: Section \ref{Sec-cpx-F} provides a brief introduction to complex Finsler geometry; Section \ref{sec-base} proves for the K\"ahler base conclusion for weakly K\"ahler--Finsler metrics;
Section \ref{sec-rig} provides the proof for Kobayashi rigidity theorems.

\section{Complex Finsler geometry}\label{Sec-cpx-F}
	
In this section, we briefly review some basic notions in complex Finsler geometry. 	
\begin{definition}[\cite{AP}]
Let $M$ be a complex manifold, a strongly pseudoconvex complex Finsler metric on $M$ is an upper semi-continuous function $F:T^{1,0}M\rightarrow [0,+\infty)$ such that

 (1) $G:=F^2$ is smooth on 	
 $\widetilde{M}:=T^{1,0}M\setminus\{\mbox{zero section}\}$;

 (2) $F(z;\lambda v)=\vert\lambda\vert F(z;v)$ for any
 $\lambda\in\mathbb{C}$ and $(z;v)\in T^{1,0}M$.
	
(3) the matrix $\left(G_{i\overline{j}}\right):=\left(\frac{\partial^2G}{\partial
v^i\partial\overline{v^j}}\right)$ is positive definite at every $(z;v)\in \widetilde{M}$.
\end{definition}

\begin{remark}
By the homogeneous of a strongly pseudoconvex complex Finsler metric $F$, $G$ is smooth on the entire bundle $T^{1,0}M$ if and only if $G$ is a smooth Hermitian quadratic metric (cf. \cite{AP}).
\end{remark}

We always use the Einstein convention in this paper.
And we adopt the convention of using a subscript to indicate differentiation of $G$ and a semicolon to distinguish derivatives with respect to the base coordinates from those with respect to the fiber coordinates. For instance, $G_{k\bar{j};i} := \dfrac{\partial^3 G}{\partial z^i \partial \bar{v}^j\partial v^k}$.
And we denote
\[\partial_i:=\frac{\partial}{\partial z^i},\,
 \dot\partial_i:=\frac{\partial}{\partial v^i},
 \partial_{\bar i}:=\frac{\partial}{\partial \bar z^i},\,
 \dot\partial_{\bar i}:=\frac{\partial}{\partial \bar v^i}.\]

In local coordinates, the connection form of the Chern--Finsler connection  $\nabla$ of $F$ is defined by:
\[\omega^i_j=G^{\bar c i}\bar\partial G_{j\bar c}
=:\Gamma^i_{j;k}dz^k+\Gamma^i_{jk}\psi^k,\]
where
\begin{equation}
\Gamma_{j;k}^{i} = \dot\partial_j \Gamma_{;k}^{i},\,
\Gamma_{;k}^{i} =  G^{\bar{s}i} G_{\bar{s};k},\,
\Gamma_{jk}^{i} = G^{\bar{s}i} G_{j\bar{s}k}.\label{con-coe}
\end{equation}
And $\Gamma_{j;k}^{i}, \Gamma_{;k}^{i}$ and $\Gamma_{jk}^{i}$
are called the horizontal, nonlinear and vertical connextion coefficients of $\nabla$, respectively.
Denoting
\[\psi^i=dv^i+\Gamma^i_{;j}dz^j,\quad \delta_i=\partial_i-\Gamma^j_{;i}\dot\partial_j, \]
then $\{dz^i;\psi^i\}$ is the locally dual coframe of $\{\delta_i;\frac{\partial}{\partial v^j}\}.$
 And \(\{\delta_i,i=1,\cdots,n\}\) forms a local frame of the horizontal bundle $\mathcal H$, a sub-bundle of tangent bundle over $\widetilde M$:
\[T^{1,0}\widetilde M=\mathcal H\oplus\mathcal V.\]

 Since $G(z;v)$ is (1,1)-homogeneous with respect to $v$ and $\bar v$, by the Euler homogeneous identity, it is easy to check (cf. \cite{AP})
\begin{eqnarray}
G(z;v)&=&G_{a}(z;v)v^a=G_{\bar a}(z;v)\bar v^a;\\
0&=&G_{ab}(z;v)v^a=G_{ab\bar c }(z;v)v^a; \\ \Gamma^i_{j;k}(z;v)v^j&=&\Gamma^i_{;k}(z;v),\,
\Gamma^i_{jk}(z;v)v^j=0.\label{AP-homo}
\end{eqnarray}

\section{K\"ahler base of weakly K\"ahler--Finsler metric}\label{sec-base}

Now we give the proof of Theorem \ref{thm}.
\begin{theorem}\label{pf}
 Every K\"ahler--Finsler metric has a K\"ahler base.
\end{theorem}
\begin{proof}
Notice that the weakly K\"ahler identity can be written as
\begin{equation}
 \delta_kG_i=\delta_iG_k. \label{condition}
 \end{equation}
For every complex Finsler metric, by lemma 2.3.4 in \cite{AP},
\begin{equation}
 \delta_kG_{\bar j}  =0,\quad  \delta_kG=0.
 \label{delta-G}
 \end{equation}
Consider the horizontal \((1,1)\)-form
\[ \beta = \sqrt{-1}\,\frac{G_iG_{\bar j}}{G}\, dz^i\wedge d\bar z^j. \tag{3} \]
It is (0,0)-homogeneous with respect to $v$, so \(\beta\) descends to \(P(T^{1,0}M)\). Although it has rank one at each fiber direction, its average will be positive definite.
Equations \eqref{condition}-\eqref{delta-G} give
\[ \delta_k\!\left(\frac{G_iG_{\bar j}}{G}\right) -\delta_i\!\left(\frac{G_kG_{\bar j}}{G}\right)
= \frac{G_{\bar j}}{G} \bigl(\delta_kG_i-\delta_iG_k\bigr)=0.  \]
Consequently \(\partial_H\beta=0\), where we denote
\[d=d_H+d_V,\, d_H=\partial_H+\bar\partial_H,\, d_V=\partial_V+\bar\partial_V,\, \partial_H \partial_H:=\delta_k\otimes dz^k,\,\partial_V=\dot\partial_k\otimes\psi^k,\]
 conjugation provides \(\bar\partial_H\beta=0\). Thus
\begin{equation}
 d_H\beta=0. \label{dH}
\end{equation}
Let \(\pi:P(T^{1,0}M)\to M\) be the canonical projection, and set \(m=n-1\). Following the methods in \cite{CSZ} by Chen--Shen--Zhao, use the globally defined closed form (cf. \cite{Kob})
\[ L=\sqrt{-1}\,\partial\bar\partial\log G =\omega_V-\Theta, \]
where
\[\omega_V=\sqrt{-1}(\log G)_{i\bar j}(z;v)dv^i\wedge d\bar v^j,\quad (\log G)_{i\bar j}=\dot\partial_i\dot\partial_{\bar j}G,\]
is vertical and positive on every projective fiber, while
\[\Theta
=\frac{\sqrt{-1}}{G}K_{i\bar j}dz^1\wedge d\bar z^j,\quad
K_{i\bar j}=-\partial_i\partial_{\bar j}G+G^{k\bar m}(\partial_i\dot\partial_{\bar m} G)(\partial_{\bar j} \dot\partial_{k} G),
\]
 is horizontal. Their fiber-integration construction and degree argument apply here.
Define
\[ \Omega=\pi_*\!\left(\beta\wedge\frac{L^m}{m!}\right) =\pi_*\!\left(\beta\wedge\frac{\omega_V^m}{m!}\right). \]
Equivalently,
\(\Omega=\sqrt{-1}\,h_{i\bar j}\,dz^i\wedge d\bar z^j\),
 with
\[  h_{i\bar j}(z) = \int_{P(T_z^{1,0}M)} \frac{G_i(z,v)G_{\bar j}(z,v)}{G(z,v)} \,\frac{\omega_V^{n-1}}{(n-1)!}.
\]
Because fiber integration commutes with \(d\), \(dL=0\), and by \eqref{dH},
\[ d\Omega = \pi_*\!\left(d_V\beta\wedge\frac{L^m}{m!}\right). \]
Every term in \(d_V\beta\wedge(\omega_V-\Theta)^m\) has vertical bidegree
\[ (k+1,k)\quad\text{or}\quad(k,k+1), \]
never \((m,m)\). Hence every term has zero integration, thus
\[ d\Omega=0. \]
Finally,
\[ h_z(X,\bar X) = \int_{P(T_z^{1,0}M)} \frac{\vert G_i(z,v)X^i\vert^2}{G(z,v)} \,\frac{\omega_V^{n-1}}{(n-1)!}\geq0,\quad 0\neq X\in T_z^{1,0}M. \]
The integrand is nonnegative. Notice at the direction \([X]\), Euler's identity for homogeneous maps gives
\[ G_i(z,X)X^i=G(z,X)>0, \quad  0\neq X\in T_z^{1,0}M. \]
It is therefore positive on an open neighborhood of \([X]\), so \[h_z(X,\bar X)>0,\quad  0\neq X\in T_z^{1,0}M.\]
 Compactness of the fibers also ensures that \(h\) is smooth.
Thus \(\Omega\) is a smooth, closed, positive \((1,1)\)-form, and \(M\) is K\"ahler.
\end{proof}

\section{ Rigidity of the K\"ahler--Finsler-Kobayashi metric }\label{sec-rig}

Convex domains are intimately connected with the pluricomplex Green function, pluripotential theory, the complex Monge--Amp\`ere equation, extremal holomorphic discs (also known as holomorphic geodesic discs) and the Kobayashi metric in the study of Lempert theory \cite{Lem81,Lem82,Lem84,RWK13,Pat,PS20}. Recently, Gaussier-Zimmer, and Dong--Wang--Wong obtained some rigidity theorems using the Kobayashi/Carath\'eodory metric for bounded strongly pseudoconvex domains using K\"ahler condition as introduced in the first section of this paper \cite{Gau-Zim,DWW}. However, the Kobayashi/Carath\'eodory metric is Finsler in general, so it's nature to consider under complex Finsler conditions, for example, using K\"ahler--Berwald condition \cite{GZ}.
In this section, we consider under the K\"ahler--Finsler condition for the Kobayashi metric of balanced pseudoconvex domains.

Here a bounded domain $D$ is called strongly convex if the local defining function of $D$ has positive definite real Hessian;
strongly pseudoconvex if the local defining function of $D$ has positive definite complex Hessian.
And a domain $D$ is called balanced if
\[D=\{\lambda z: z\in D, \lambda \in \mathbb D\}.\]
Balancedness identifies the domain with its Kobayashi indicatrix at the origin. Indeed, if $\mu_D$ denotes the Minkowski functional of a balanced bounded pseudoconvex domain $D$, then (by \cite{Barth,Suzuki}, they called balanced domain as complete circular domain)
\[
 F_K(0;v)=\mu_D(v),\quad
 D=\{v\in\mathbb C^n:\mu_D(v)<1\}.
\]
 Consequently, the function
$g(v):=F_K(0;v)^2=\mu_D(v)^2$ is smooth  and plurisubharmonic on $T^{1,0}_0D\setminus \{0\}$. On the other hand, we can weaker the boundary condition, with corresponding assumption for the metric $F_K^2$ near the origin, since the boundary of the domain is defined by the norm $g$.
Hence the the convexity and regularity of the metric and the domain contral each other.

In the first subsection, we treat the strongly convex case,
extending the arguments for the special examples considered
in the author's earlier work \cite{G}. This approach uses
classical Lempert theory, particularly stationary complex
geodesic discs.
In the next subsection, we consider for balanced, smooth,
bounded, strongly pseudoconvex domains. Where the propagation theorem of Patrizio and Spiro \cite{PS20} provides smoothly varying stationary infinitesimal extremal discs and smoothness of the Kobayashi metric away from the zero section over a neighborhood
of the center.

\subsection{The strongly convex case}\label{Lem-conv}

Recall that, for a complex manifold $M$, and \((z,v) \in \widetilde{M}\), a holomorphic map \(f_{z, v}: \mathbb D \to M\) is called extremal at \((z, v)\), if
\[f_{z,v}(0)=z,\,\dot f_{z,v}(0)=\lambda v,\]
with \(\vert \lambda\vert\) being maximal (in case $\Omega$ is balanced, we can assume $\lambda\in \mathbb R$).
In \cite{Lem81},  Lempert obtained the following properties for smooth, bounded, strongly convex domains (cf. also \cite{Lem82,Lem84}, and \cite{RWK13} for generalization results about bounded convex domain without boundary regularity assumption), for any $(z,v)\in \widetilde \Omega$:

(L1) there is a unique extremal map $f_{z,v}:\mathbb D\to \Omega$ such that $ f(0)=z,\, \dot f(0)=\lambda v$.

(L2) $f_{z,v}$ is the unique geodesic complex curve through $(z,v)$, i.e., it is an isometry with respect to the Poincar\'e distance of $\mathbb D$ and the Kobayashi distance of $\Omega$, thus  $F_K(z;v)=\frac1\lambda$.

(L3) Smooth dependence: $(z,v)\mapsto f_{z,v}$ and
$(z,v)\mapsto\lambda(z,v)$ smoothly depend on $\widetilde \Omega.$ And $f_{z,v}$ extends smoothly to $\bar{\mathbb D}$, embedding $\bar{\mathbb D}$ into $\bar \Omega$, $f(\partial\mathbb D)\subset \partial\Omega$.

(L4) Dual-map criterion (stationary geodesic disk): a holomorphic immersion $\varphi:\mathbb D\to \Omega$, smooth up to the boundary, with $\varphi(\partial\mathbb D)\subset \partial \Omega$, is extremal if and only if there exists a real map $p>0$ defined (for almost every point) on $\partial\mathbb D$, such that
\[H_j(\zeta):=\zeta p(\zeta)\rho_{;j}(\varphi(\zeta)),\,j=1,2,\, \rho_{;j}=\partial_j\rho,\]
extends holomorphically to $\mathbb D$, where $\rho$ is the defining function of $\Omega$ ( for our case, $\rho$=g-1, thus $\rho_{;j}=g_j$).

And we also need the following condition for K\"ahler--Finsler metric.
\begin{theorem}[\cite{G}]\label{G}
A strongly pseudoconvex weakly K\"ahler--Finsler metric is K\"ahler--Finsler if and only if
$ \mathbb Q_{i\bar j}\equiv0,$
where $\mathbb Q_{i\bar l}$ is induced from the complex geodesic spray coefficients $\mathbb G^i$:
\[
\mathbb G^i:=\Gamma^i_{;j}v^j,\quad
\mathbb G^i_{\bar l}:=\dot\partial_{\bar l}\mathbb G^i,\quad
\mathbb Q_{j\bar l}:=G_{ij}\mathbb G^i_{\bar l},
\]
\end{theorem}

In the following we use notions $\dot \psi:=\frac{\partial \psi}{\partial\zeta}$ for $\zeta\in \mathbb D,$ $\lambda':=\frac{\partial \lambda}{\partial\varepsilon}$ for $\varepsilon\in \mathbb R,$ and we use lower indices to represent components of points and vectors for calculation convenience.
Now we can establish the following rigidity theorem.

\begin{theorem}\label{rigi-KF}
For every balanced, smooth, bounded, strongly convex domain
$D\subset\mathbb C^n$, its Kobayashi metric is K\"ahler--Finsler
if and only if $D$ is a Hermitian ellipsoid, linearly biholomorphic to the unit ball $\mathbb B^n$.
\end{theorem}
\begin{proof}
At first, we consider for the necessary direction.  Let
\[G(z;v)=F_K(z;v)^2,\quad g(v)=G(0;v).\]
Since $D$ is bounded, balanced and convex, the linear discs
$\varphi_0(\zeta)=\zeta u$, $u\in\partial D$, are the complex geodesics through the origin, and the Kobayashi indicatrix at $0$ coincides with $D$:
\[D=\{z\in \mathbb C^n:g(z)<1\}\]
(cf. Proposition 3.5.3 and Remark 11.1.3(b) in \cite{Jar-Pfl}); in
particular $\rho:=g-1$ is a smooth defining function of $D$ near
$\partial D$.
By Lempert's theorem (\cite{Lem81}), $G$ is smooth on
$\widetilde D:=T^{1,0}D\setminus\{\text{zero section}\}$. By condition, $F_K$ is K\"ahler--Finsler,
hence $(g_{i\bar j})$ is positive
definite on $\mathbb C^n\setminus\{0\}$, and the $(1,1)$-homogeneity of $g$ gives Euler's identities
\begin{equation}
g_i(v) v^i=g,\quad g_{ij}(v)v^j=0,\quad g_{i\bar j}(v)\bar v^j=g_{i}(v).
\label{euler}
\end{equation}

\textbf{Step 1.}
 The central stationary disc is $\varphi_0(\zeta)=\zeta u$ for any fixed $u\in \partial D$, i.e., $g(u)=1$ (cf. Proposition 3.11, Proposition 3.5.3 and Remark 11.1.3(b) in \cite{Jar-Pfl}).
Here we provide another check for the extremal map at the origin, for later use.

For \(g(u)=1\), Euler formula for homogeneous maps gives,
   \[   g_i(\zeta u)=\bar\zeta g_i(u),\quad
   g_{ij}(\zeta u)=\bar\zeta^2g_{ij}(u),\quad
   g_{i\bar j}(\zeta u)=g_{i\bar j}(u),\quad \vert\zeta\vert=1.
   \]
   Thus the lift of \(\varphi_0(\zeta)=\zeta u\), together with
the real valued map $p(\zeta):=\frac{1}{\vert\zeta\vert^2 }$ makes
\[
 \zeta p(\zeta) g_i(\zeta u)
=\zeta\frac{1}{\vert\zeta\vert^2 } \bar\zeta g_i( u)
=g_i( u)
  \]
a constant map with respect to $\zeta$ near the boundary $\partial\mathbb D$, thus extends holomorphically.
This proves central stationarity directly using (L4). And
$\varphi_0(\zeta)=\zeta u$ is the holomorphic geodesic disc for
$(0,u)$, $u\in\partial D$. Below we consider with fixed $u$.

\textbf{Step 2.} We consider the first variation of $\varphi_0$.
Let $\varepsilon\in\mathbb R$, and denote $\varphi_\varepsilon$  the normalized holomorphic geodesic disks $\varphi_{\varepsilon w,\lambda(\varepsilon)u}$ satisfying
\[
\varphi_\varepsilon(0)=\varepsilon w,\quad
\varphi_\varepsilon'(0)=\lambda(\varepsilon)u,\quad
\lambda(\varepsilon)>0.
\]
By (L3), the smooth dependence of normalized complex geodesics, the  family $\varphi_\varepsilon$ depends smoothly on \(\varepsilon\).
By (L2), $\varphi_\varepsilon$ is extremal at $(\varepsilon w,u)$, and by definition,
\begin{equation}
G(\varepsilon w;u)=\lambda(\varepsilon)^{-2}.\label{Glambda}
\end{equation}
For the varying extremal discs, there exist positive $p_\varepsilon$ on $\partial \mathbb D$, such that the pairing
$\widetilde\varphi_{\varepsilon,i}
:=\zeta\,p_\varepsilon(\zeta)g_i(\varphi_\varepsilon(\zeta))$ extends holomorphically to $\mathbb D$ by (L4).

Put
\[
A_\varepsilon(\zeta)
=\zeta g_i(\varphi_\varepsilon(\zeta))
       \varphi_\varepsilon'{}^i(\zeta).
\]
Because every \(\varphi_\varepsilon\) maps \(\partial\mathbb D\) into \(\partial D=\{g=1\}\), we have
\[
g(\varphi_\varepsilon(\zeta))=1,\quad \vert\zeta\vert=1.
\]
Differentiating the boundary equation $g(\varphi_\varepsilon(e^{i\theta}))\equiv1$ in \(\theta\), where \(\zeta=e^{i\theta}\), gives
\[
0=\frac{d}{d\theta}g(\varphi_\varepsilon(e^{i\theta}))
  =-2\operatorname{Im}A_\varepsilon(\zeta).
\]
Thus \(A_\varepsilon\) is real on \(\partial\mathbb D\). Consequently,
\[
\widetilde\varphi_{\varepsilon,i}
\varphi_\varepsilon'{}^i=p_\varepsilon A_\varepsilon
\]
is holomorphic with real boundary values, hence constant. Near \(\varepsilon=0\) this constant is positive.
By rescaling normalization, we can assume $p_\varepsilon A_\varepsilon\equiv1$ on $\partial \mathbb D$, i.e.,
\[
p_\varepsilon(\zeta)
=\frac{1}{\zeta\,g_i(\varphi_\varepsilon(\zeta))
\varphi_\varepsilon'{}^i(\zeta)},\quad \zeta\in\partial \mathbb D.
\]
The denominator is real, thus positive for small $\varepsilon$ since it equals $1$ when $\varepsilon=0$, hence $p_\varepsilon(\zeta)$  depends smoothly on small $\varepsilon$.

 Write
\[ \varphi_\varepsilon(\zeta)
=\zeta u+\varepsilon\psi(\zeta)+\tfrac12\varepsilon^2\chi(\zeta)
+\cdots,  \quad \lambda(0)=1,\]
so that
\[\psi(0)=w,\quad \psi'(0)=\lambda'(0)u,\quad
\chi(0)=0,\quad \chi'(0)=\lambda''(0)u.\]
Fix such \(\zeta\in\partial \mathbb D\), differentiate
\[
g(\varphi_\varepsilon(\zeta))=1,\quad \vert\zeta\vert=1
\]
 with respect to the real parameter \(\varepsilon\) at \(0\) gives
\[
0=g_i(\zeta u)\psi^i(\zeta)
  +g_{\bar i}(\zeta u)\overline{\psi^i(\zeta)}
 =2\operatorname{Re}\!\left(g_i(\zeta u)\psi^i(\zeta)\right).
\]
By homogeneous of $g_i$, $g_i(\zeta u)=\bar\zeta g_i( u)$, thus
\[   \operatorname{Re}\bigl(\bar\zeta g_i(u)\psi^i\bigr)=0.  \]
  Its constant mode implies
\[ \lambda'(0)=0,\quad g_i a_2^i=-g_{\bar l}\bar w^l,  \quad g_ia_m^i=0\quad(m\ge3),
\]
and the remaining modes give
   \[\psi=w+\sum_{m\ge2}a_m\zeta^m.\]
If $p_\varepsilon=1+\varepsilon q+\cdots$, the first variation of the lift
$\zeta p_\varepsilon g_i(\varphi_\varepsilon)$ is
\[qg_i+\bar\zeta g_{ij}\psi^j +\zeta g_{i\bar j}\overline{\psi^j}.\]
Its negative Fourier modes vanish, since the lift extends holomorphically to $\mathbb D$.
Set $a=a_2$ and write $q=\sum q_m\zeta^m$. For $m\ge2$, the
$\zeta^{-m}$-coefficient is
\[
q_{-m}g_i+g_{i\bar j}\overline{a_{m+1}^j}=0.
\]
Contracting with $u^i$, using \eqref{euler} ($u^ig_i=1$,
$u^ig_{i\bar j}=g_{\bar j}$) and
$g_{\bar j}\overline{a_{m+1}^j}=\overline{g_ja_{m+1}^j}=0$, gives
$q_{-m}=0$; the invertibility of $(g_{i\bar j})$ then gives $a_{m+1}=0$.
The $\zeta^{-1}$-mode gives
\[q_{-1}g_i+g_{ij}w^j+g_{i\bar j}\bar a^j=0;\]
contracting with $u^i$, using $u^ig_{ij}=0$ and
$g_ia^i=-g_{\bar l}\bar w^l$, yields
\[q_{-1}=g_iw^i,\quad  g_{i\bar j}\bar a^j=-g_{ij}w^j-g_i\,q_{-1},\]
and solving the second equation (conjugate and multiply by $g^{\bar ai}$,
using $g^{\bar ai}(u)g_{\bar a}(u)=u^i$) gives
\begin{equation}
a^i=-g^{\bar a i}\overline{g_{al}}\bar w^l-u^i g_{\bar l}\bar w^l
=-B^{i}_{\bar l}\bar w^l-u^i g_{\bar l}\bar w^l,\quad B^i_{\bar l}:=g^{\bar a i}\overline{g_{al}}.
\end{equation}
Therefore
\[\varphi_\varepsilon(\zeta)
=\zeta u+\varepsilon(w+a\zeta^2)+O(\varepsilon^2).\]

\textbf{Step 3.} Second boundary variation and the expansion of $G(z;v)$.
Differentiating the identity $g(\varphi_\varepsilon(\zeta))\equiv1$
($\vert\zeta\vert=1$) twice at $\varepsilon=0$ and using the homogeneity relations
of Step 1,
\[
0=\frac12\frac{d^2}{d\varepsilon^2}\Big\vert_{\varepsilon=0}
g(\varphi_\varepsilon(\zeta))
=\operatorname{Re}\bigl(\bar\zeta g_i\chi^i\bigr)
+\operatorname{Re}\bigl(\bar\zeta^2 g_{ij}\psi^i\psi^j\bigr)
+g_{i\bar j}\psi^i\bar\psi^j,\quad \vert\zeta\vert=1.
\]
Since $\chi=\lambda''(0)u\zeta+O(\zeta^2)$, its constant Fourier mode is
\begin{equation}
\lambda''(0) + 2\operatorname{Re}(g_{ij}w^ia^j)
 + g_{i\bar j}w^i\bar w^j + g_{i\bar j}a^i\bar a^j=0.\label{lambda2}
\end{equation}
Denoting
\[C_{k\bar l}:=g_{km}B^i_{\bar l}
=g_{km}g^{\bar a m}\overline{g_{al}},\]
by Euler formula for homogeneous maps as given in \eqref{euler},
one obtains
\begin{equation}
g_{ij}(u)u^j=0,\quad
 g_iB^i_{\bar l}=\bar v^kg_{i\bar k}g^{\bar ji}\overline{g_{jl}}=\bar v^j\overline{g_{jl}}=0.
\label{eulerBC}
\end{equation}
Thus
\[  g_{ij}w^ia^j=-C_{k\bar l}w^k\bar w^l, \quad
g_{i\bar j}a^i\bar a^j=C_{k\bar l}w^k\bar w^l+\vert g_kw^k\vert^2,\]
hence by \eqref{lambda2}, one has
\[-\lambda''(0)=g_{k\bar l}w^k\bar w^l-C_{k\bar l}w^k\bar w^l
+\vert g_kw^k\vert^2.\]
Balancedness gives $z\mapsto-z\in\mathrm{Aut}(D)$, so $G(-z;-v)=G(z;v)$;
together with the homogeneity about fibers $G(-z;-v)=G(-z;v)$,
one has
\[G(z;v)=G(-z;v),\]
thus $G$ is an even map. So all odd-order
Taylor terms of $G(\cdot;v)$ at $0$ vanish, and $\lambda(\varepsilon)$ is
even. Using \eqref{Glambda},
\[G(\varepsilon w;u)=\lambda(\varepsilon)^{-2}
=1+\varepsilon^2
\left[g_{k\bar l}w^k\bar w^l
-C_{k\bar l}w^k\bar w^l
+\left\vert g_kw^k\right\vert^2
\right]+O(\varepsilon^4).\]
Since this holds for every $w$, the second-order Taylor polynomial of
$G(\cdot;u)$ at $0$ is determined: a putative term
$A_{kl}z^kz^l+\overline{A_{kl}}\bar z^k\bar z^l$ would contribute
$2\varepsilon^2\operatorname{Re}(A_{kl}w^kw^l)$ to the $\varepsilon^2$-coefficient
above, forcing $\operatorname{Re}(A_{kl}w^kw^l)\equiv0$ for all $w$, hence
$A_{kl}=0$. Thus
\[G(z;u)=1+\bigl[g_{k\bar l}-C_{k\bar l}+g_kg_{\bar l}\bigr]z^k\bar z^l
+O(\vert z\vert^4).\]
Rescaling $u=\frac{v}{\sqrt{g(v)}}\in\partial D$: since $g_{k\bar l}$ and
$C_{k\bar l}$ are $0$-homogeneous while
$g_k(u)g_{\bar l}(u)=g_k(v)g_{\bar l}(v)/g(v)$, the homogeneity
$G(z;v)=g(v)G(z;u)$ yields
\begin{equation}
G(z;v)=g(v)+T_{k\bar l}(v)z^k\bar z^l+O(\vert z\vert^4),
\quad
T_{k\bar l}=g(g_{k\bar l}-C_{k\bar l})+g_k g_{\bar l}.\label{expansion}
\end{equation}
The expansion is locally uniform in $v\ne0$; since $G$ is smooth on
$\widetilde D$, after one base derivative the remainder and its required fiber derivatives are $O(\vert z\vert^3)$.
Notice, from \eqref{euler} and \eqref{eulerBC}, one has
\begin{equation}
T_{k\bar l}v^k=2g\,g_{\bar l}.\label{Tcontract}
\end{equation}

\textbf{Step 4.} Vanishing of the obstruction.
We use the complex geodesic spray and consequent notions as before:
\[
\mathbb G^i=\Gamma^i_{;j}v^j=G^{\bar s i}G_{\bar s;j}v^j,\quad
\mathbb G^i_{\bar l}=\dot\partial_{\bar l}\mathbb G^i,\quad
\mathbb Q_{j\bar l}=G_{ij}\mathbb G^i_{\bar l}.
\]
From \eqref{expansion} and \eqref{Tcontract},
\begin{eqnarray*}
G_{\bar s;j}v^j
&=&\dot\partial_{\bar s}\bigl(T_{j\bar l}v^j\bigr)\bar z^l+O(\vert z\vert^3)\\
&=&\dot\partial_{\bar s}(2g\,g_{\bar l})\,\bar z^l+O(\vert z\vert^3)\\
&=&2\bigl(g_{\bar s}g_{\bar l}+g\,g_{\bar s\bar l}\bigr)\bar z^l
+O(\vert z\vert^3),
\end{eqnarray*}
using $G^{\bar s i}(z;v)=g^{\bar s i}+O(\vert z\vert^2)$,
$g^{\bar s i}g_{\bar s}=v^i$, and $g^{\bar s i}g_{\bar s\bar l}=B^i_{\bar l}$,
\begin{equation}
\mathbb G^i(z;v)
=2\left(v^i g_{\bar l}+gB^i_{\bar l}\right)\bar z^l
+O(\vert z\vert^3).\label{spray}
\end{equation}
By $G_{ij}(z;v)=g_{ij}(v)+O(\vert z\vert^2),$ $g_{ij}v^i=0$, and \eqref{spray}, one has
\begin{equation}
 \mathbb Q_{j\bar k}(z;v)
 =G_{ij}(z;v)\mathbb G^i_{\bar k}(z;v)
 =2\left(g_{\bar k}C_{j\bar l}
   +g\,g_{ij}\dot\partial_{\bar k}B^i_{\bar l}
   \right)\bar z^l+O(\vert z\vert^3).\label{KF-Obstruction}
 \end{equation}
Notice
$B^i_{\bar l}\bar v^l=g^{\bar ji}\overline{g_{jl}}\bar v^l=0,$
one has
\begin{equation}
  \bigl(\dot\partial_{\bar k}B^i_{\bar l}\bigr)\bar v^l
=\dot\partial_{\bar k}\bigl(B^i_{\bar l}\bar v^l\bigr)-B^i_{\bar k}
=-B^i_{\bar k}.\label{diff-eq}
\end{equation}
Evaluation at $z=\varepsilon v$, together with $C_{j\bar l}\bar v^l=g_{km}g^{\bar a m}\overline{g_{al}v^l}=0$, and substituting \eqref{diff-eq} to \eqref{KF-Obstruction},
provides
\begin{equation}
\mathbb Q_{j\bar l}(\varepsilon v;v)
=-2\bar\varepsilon\,g(v)C_{j\bar l}(v)
+O(\vert\varepsilon\vert^3).\label{expansion-Q}
\end{equation}

By Lempert's theory and Corollary 3.4 of Abate--Patrizio \cite{AP2}, $F_K$ is a strongly convex weakly K\"ahler--Finsler metric.
If $F_K$ is K\"ahler--Finsler, Theorem \ref{G} gives
\begin{equation}
0=\mathbb Q_{j\bar l}(\varepsilon v;v)
=-2\bar\varepsilon\,g(v)C_{j\bar l}(v)
+O(\vert\varepsilon\vert^3).\label{KF-eq}
\end{equation}
Notice the quantity $g(v)C_{j\bar l}(v)$ in \eqref{KF-eq} is independent of small $\varepsilon$ and $g(v)>0$, \eqref{KF-eq} provides
\[C_{j\bar l}(v)=0, \quad \forall v\neq0.\]
Notice $ C_{k\bar l}=g_{km}g^{\bar a m}\overline{g_{al}},$
and  $(g^{\bar a m})=(g_{i\bar j})^{-1}$ is a positive definite matrix, thus
\[ (g_{ij})(v)=0, \quad \forall v\neq0.\]
Since fiber derivatives commute,
$g_{i\bar jk}=\dot\partial_{\bar j}(g_{ik})=0$, and similarly all third and higher fiber derivatives of $g$ vanish; hence $(g_{i\bar j})$ is a constant positive Hermitian matrix, and Euler's identity gives
\[g(v)=g_{i\bar j}v^i\bar v^j.\]
Then $g$ is Hermitian quadratic, and $D=\{g<1\}$ is a Hermitian ellipsoid. And the strongly pseudoconvexity makes the Hermitian ellipsoid linearly biholomorphic to the ball.

For the sufficient direction: if $D$ is biholomorphic to the unit ball $\mathbb B^n$, by the invariance of the Kobayashi metric under biholomorphisms, the Kobayashi metric of $D$ is K\"ahler since the Kobayashi metric $F_K$ of $\mathbb B^n$ is K\"ahler.
\end{proof}

\subsection {The pseudoconvex case}
By the proof in Subsection \ref{Lem-conv}, the K\"ahler--Finsler rigidity results extends to any bounded Lempert domains (and Lempert manifolds) (see \cite{ALY19,PS20} for their definions) that have stationary geodesic disks (satisfying (L4)). 
Here we consider  generally case, on strongly pseudoconvex and even pseudoconvex (formally) domains, for the Kobayashi rigidity with K\"ahler--Finsler condition.

In \cite{Pang}, M.-Y. Pang proved that, any balanced, smooth, bounded, strongly pseudoconvex domain $D$ admits an open subset $U\subset D$, such that each $y\in U$ is the center of an smooth Monge--Amp\`ere exhaustion $\tau^y$. And the family of exhaustions $\tau^{(x)}$, considered as collections of maps parameterized by the uniquely associated centers   $x\in \mathbb U$, depends smoothly  on such parameter.
Thus the Kobaysashi metric of $D$ restricted on each center $x\in U$ of $\tau^{(x)}$, in local holomorphic coordinates, is given by (\cite{Pat87})
\[\displaystyle \kappa_D(x;v)
=\lim_{r\downarrow0}\frac{\sqrt{\tau^{(x)}(x+rv)}}{r},
\quad v\ne0.\]
thus $F_K$ depends smoothly on $x\in U$ \cite{Pang}. And $F_K$ is
 strong pseudoconvexity because of the foliation theory also gives strong pseudoconvexity of
 the indicatrices
\[
I_x=\{v:\kappa_D(x;v)<1\}.
\]
In \cite{PS20}, Patrizio and Spiro generalized Pang's results as follows.

\begin{theorem} [\cite{PS20}, Theorem 1.1 and Theorem 2.9]\label{PS-MA}
Let $ D \subset \mathbb C^n$ be a  smoothly bounded strongly pseudoconvex domain,  admitting a  Monge--Amp\`ere exhaustion
$\tau: D \to [0,1)$,  centered at $x_o \in D$, which  is smoothly  extendible up  to the boundary. Then
there exists  an open neighborhood $\mathbb U\subset D$ of $x_o$, such that for each  $x \in \mathbb U$ there exists a  Monge--Amp\`ere exhaustion $\tau^{(x)}: M \to [0,1)$,  smoothly extendible up to the boundary and centered at $x$.
The family of exhaustions \(\{\tau^{(x)}\}_{x\in\mathbb U}\) depends smoothly on the center \(x\).
\end{theorem}

The preliminary lemma below provides these properties in the form needed for the first and second variations in the proof of Theorem \ref{rigi-KF} in last subsection.
Which is a consequence of \cite[Theorem~2.9 and Section~2.2]{PS20}, together with the identity $F_K(0;\cdot)=\mu_D$, as a special exhaustion, for balanced pseudoconvex domains.
We include a proof to make the parameter dependence and the normalization of the holomorphic lifts explicit.

\begin{lemma}\label{lem:balanced}
Let $D\subset\mathbb C^n$, $n\geq2$, be a balanced, bounded pseudoconvex domain.
Suppose that $g(v):=F_K(0;v)^2$ is smooth and strongly pseudoconvex on $\mathbb C^n\setminus\{0\}$. Then $D=\{g<1\}$ has smooth strongly pseudoconvex boundary.
There is a neighborhood $U$ of $0$ such that one can choose
discs $f_{z,v}\in\mathcal O(\mathbb D, D)\cap C^\infty(\overline{\mathbb D},\overline D)$,
depending smoothly on $(z,v)\in U\times(\mathbb C^n\setminus\{0\})$, with
\[
 f_{z,v}(0)=z,\quad
 f_{z,v}'(0)=\frac{v}{F_K(z;v)},\quad
 f_{z,v}(\partial\mathbb D)\subset\partial D.
\]
They admit positive smooth boundary multipliers $p_{z,v}$ such that the functions
\[
 H_{z,v,i}(\zeta)
 :=\zeta p_{z,v}(\zeta)g_i(f_{z,v}(\zeta)),\quad \zeta\in\partial\mathbb D,\,i=1,\ldots,n,
\]
extend holomorphically to $\mathbb D$ and satisfy
\[
 H_{z,v,i}f_{z,v}'{}^i\equiv1.
\]
At the balanced center,
\[
 f_{0,v}(\zeta)=\frac{\zeta v}{\mu_D(v)}.
\]
\end{lemma}
\begin{proof}
Let $\mu_D$ denote the Minkowski functional of $D$. Since $D$ is
balanced and pseudoconvex, the restriction of $F_K$ to the tangent space at the origin coincides with $\mu_D$:
\[
g(v):=F_K^2(0;v)=\mu_D^2(v).
\]
Thus, \[D=\{g<1\},\quad \partial D=\{g=1\}.\]
By the Euler formula for homogeneous function $g$, $g_{i\bar j}(b)b^i\bar b^j=g(b)=1,\, \forall b\in \partial D,$ thus $dg\neq 0$ on $\partial D$. Moreover, $(g_{i\bar j})$ is positive definite. Hence $\rho:=g-1$ is a smooth strictly plurisubharmonic defining function near $\partial D$, so $D$ has a smooth and strongly pseudoconvex boundary.

The squared Minkowski functional of a balanced, smooth, bounded, strongly pseudoconvex domain is a smooth Monge--Amp\`ere exhaustion,
with smoothness at the center understood on the blow-up
(see Section 2.2 in \cite{PS20}). Thus $g\vert_{D}$ is such an
exhaustion centered at $0$, smoothly extendible up to $\partial D.$

By Theorem \ref{PS-MA}, there exists a neighborhood $U$ of $0$ and Monge--Amp\`ere exhaustions $\tau_z$ centered at $z\in U$,
with $\tau_0=g$, depending jointly smoothly on the center and the
spatial variable away from the moving center, including
with this regularity extending up to $\partial D$. Each $\tau_z$ pulls back smoothly to the blow-up at $z$. The equality $\tau_0=g$ follows from the uniqueness of a Monge--Amp\`ere exhaustion with a prescribed center.

Applying properties (i)--(iii) of the circular representation
associated with $\tau_z$
\cite[Section 2.2, p. 1299]{PS20}, and parametrizing its radial discs by their boundary value at $\zeta=1$, we obtain, for every
$z\in U$ and $b\in\partial D$, a mapping
\[\varphi_{z,b}\in\mathcal O(\mathbb D,D)\cap
C^\infty(\overline{\mathbb D},\overline D)\]
 such that
\[
\varphi_{z,b}(0)=z,\quad \varphi_{z,b}(1)=b,\quad
\tau_z(\varphi_{z,b}(\zeta))=\vert\zeta\vert^2.
\]
On the punctured disc, the complex tangent line to its image coincides with the kernel of $\partial\bar\partial\log\tau_z$.
And these leaves are infinitesimal Kobayashi extremals at their centers by the equation (2.3) in \cite{PS20}, that is,
\[
F_K\bigl(z;\varphi_{z,b}'(0)\bigr)=1.
\]
At $z=0$, the kernel computed above is radial, hence the homogeneity of $g$ implies that the kernel of $\partial \bar\partial \log g$ is radial. Hence
\[
\varphi_{0,b}(\zeta)=\zeta b.
\]
We next justify the joint smoothness of \((z,b,\zeta)\mapsto\varphi_{z,b}(\zeta)\) for \(\zeta\in\overline{\mathbb D}\). On a fixed collar neighborhood of \(\partial D\), define the \((1,0)\) vector field \(Z_z\) by
\[
(\tau_z)_{i\bar j}Z_z^i=(\tau_z)_{\bar j}.
\]
After shrinking \(U\), these vector fields are jointly smooth. 
Since \(\tau_0=g\) has positive-definite complex Hessian near \(\partial D\), shrinking \(U\) if necessary ensures that these vector fields are well defined and jointly smooth in \(z\) and the spatial variable. Near \(\partial\mathbb D\), the leaf identities imply
\[
\zeta\varphi_{z,b}'(\zeta)
=Z_z(\varphi_{z,b}(\zeta)),
\quad
\varphi_{z,b}(e^{i\theta})
=\operatorname{Fl}_{\,iZ_z-i\overline{Z_z}}^\theta(b).
\]
Thus the boundary loops are jointly smooth in $(z,b,\theta)$.
Since each component of $\varphi_{z,b}$ is the Poisson extension
of its boundary values, the regularity of the Poisson extension
operator implies that $(z,b,\zeta)\mapsto\varphi_{z,b}(\zeta)$ is jointly smooth on $\overline{\mathbb D}$.
In particular, $(z,b)\mapsto\varphi_{z,b}'(0)$ is smooth.

We now select the prescribed tangent direction. Set
\[
A_z(b)=\varphi_{z,b}'(0),\quad
E_z(b)=\frac{A_z(b)}{\|A_z(b)\|}.
\]
Since
\[
A_0(b)=b,\quad E_0(b)=\frac b{\|b\|},
\]
the map $E_0:\partial D\to S^{2n-1}$ is a diffeomorphism.
By compactness, after shrinking $U$ once more, $A_z$ is nowhere
zero and every $E_z$ is a local diffeomorphism. Each $E_z$ is
proper and hence is a covering of $S^{2n-1}$. Since $\partial D$
is connected and $S^{2n-1}$ is simply connected for $n\geq2$,
each $E_z$ is a diffeomorphism. The parameter-dependent inverse
function theorem gives smooth dependence of its inverse $E_z^{-1}$ on $z$ and its argument.

For $v\neq0$, define
\[
b(z,v)=E_z^{-1}\left(\frac v{\|v\|}\right),
\quad
f_{z,v}=\varphi_{z,b(z,v)}.
\]
These discs depend jointly smoothly on $(z,v,\zeta)$ up to
$\overline{\mathbb D}$, and
\[
f_{z,v}(0)=z,\quad
f_{z,v}(\partial\mathbb D)\subset\partial D.
\]
By construction, $f_{z,v}'(0)=a(z,v)v$ for some $a(z,v)>0$.
Center extremality and homogeneity of $F_K$ give
\[
1=F_K\bigl(z;f_{z,v}'(0)\bigr)
 =a(z,v)F_K(z;v).
\]
Hence
\[
f_{z,v}'(0)=\frac{v}{F_K(z;v)}.
\]
At $z=0$, the selected boundary point is
$b(0,v)=\frac{v}{\mu_D(v)}$, so
\[
f_{0,v}(\zeta)=\frac{\zeta v}{\mu_D(v)}.
\]

It remains to construct the normalized holomorphic lifts.
Fix $(z,v)$, write $f=f_{z,v}$, and put $u_z=\log\tau_z$.
The Monge--Amp\`ere kernel condition gives
\[
\frac{\partial}{\partial\bar\zeta}
 \bigl((u_z)_{;i}(f(\zeta))\bigr)
=(u_z)_{;i\bar j}(f(\zeta))
 \overline{f'{}^j(\zeta)}
=0
\]
on $\mathbb D\setminus\{0\}$. Therefore
\[
H_i(\zeta)=\zeta(u_z)_{;i}(f(\zeta))
\]
is holomorphic on the punctured disc.

The estimates near the center for a smooth Monge--Amp\`ere exhaustion are
\[
\tau_z(w)\asymp\|w-z\|^2,\quad
\partial\tau_z(w)=O(\|w-z\|).
\]
To justify the second estimate, observe that smoothness on the blow-up and quadratic vanishing imply, in blow-up coordinates, first derivatives of order $O(\vert\lambda\vert)$ in the radial variable and $O(\vert\lambda\vert^2)$ in the angular variables; transforming back to the original coordinates gives the asserted estimate.
Since $\tau_z(f(\zeta))=\vert\zeta\vert^2$, we have
$\|f(\zeta)-z\|\asymp\vert\zeta\vert$, and hence
\[
\vert H_i(\zeta)\vert
=\left\vert
 \frac{\zeta(\tau_z)_{:i}(f(\zeta))}{\vert\zeta\vert^2}
 \right\vert
=O(1).
\]
Thus the singularity at $0$ is removable.

On $\partial D$, both $\tau_z$ and $g$ equal $1$, so
\[
d\tau_z=c_z\,dg, \quad c_z(b)=\frac12\,d\tau_z\vert_b(b).
\]
Here we used $dg_b(b)=2$ by Euler formula for homogeneous mapping.
Here the proportionality follows because both differentials annihilate \(T_b\partial D\), and the formula for \(c_z\) follows from Euler's formula \(dg_b(b)=2\). The function \(c_z(b)\) is jointly smooth and $c_0\equiv1$; hence, after shrinking $U$, it is
positive. Define
\[
p_{z,v}(\zeta)=c_z(f_{z,v}(\zeta)),
\quad \zeta\in\partial\mathbb D.
\]
These multipliers are positive and jointly smooth in
$(z,v,\zeta)$. Since $\tau_z=1$ on the boundary $\partial D$,
 the boundary values of the lifts are
\[
H_i(\zeta)
=\zeta(\tau_z)_{:i}(f_{z,v}(\zeta))
=\zeta p_{z,v}(\zeta)g_i(f_{z,v}(\zeta)).
\]
Finally, for $\zeta\neq0$,
\[
H_i(\zeta)f_{z,v}'{}^i(\zeta)
=\zeta\,\frac{\partial}{\partial_\zeta}
 [\log\tau_z(f_{z,v}(\zeta))]
=\zeta\,\frac{\partial}{\partial_\zeta}[\log\vert\zeta\vert^2]
=1.
\]
Both sides extend holomorphically across $0$, so the identity
holds throughout $\mathbb D$.
\end{proof}

Together Lemma \ref{lem:balanced} with the variation arguments for Theorem \ref{rigi-KF}, we obtain the following rigidity result.
\begin{theorem}\label{thm:balanced-rigidity}
Let $D\subset\mathbb C^n$, $n\geq2$, be a balanced bounded pseudoconvex domain. Then its Kobayashi metric is a smooth, strongly pseudoconvex, K\"ahler--Finsler metric restricted on a neighborhood of the origin
if and only if $D$ is a Hermitian ellipsoid, linearly biholomorphic to $\mathbb B^n$.
\end{theorem}

At the end, we remark that, weakly K\"ahler--Finsler condition makes no such rigidity for the Kobayashi metric of bounded domain, since the Kobayashi metrics of smooth bounded strongly convex domains are always weakly K\"ahler--Finsler.


As a corollary, for every balanced, smooth, bounded, strictly convex, and strongly pseudoconvex domain which does not biholomorphic to the unit ball, its Kobayashi metric is a non-K\"ahler--Finsler, weakly K\"ahler--Finsler metric.

\vskip0.4cm
\noindent
\textbf{Acknowledgement.}\ {\small
The author would like to thank Prof. Xiaoshan Li for his lectures on the Cheng--Yau conjecture at AMSS CAM in 2026, and Prof. W. Zwonek for his lecture on the extension of Lempert's theorem at CNU last year. Gratitude is also extended to Prof. Xiangyu Zhou and Prof. Liyou Zhang for inviting the two aforementioned speakers respectively.}

\vskip0.4cm
\noindent
\textbf{AI Use Disclosure.}\ {\small
The research direction and initial mathematical framework were developed by the author. Generative AI (GPT-6.1 Sol, OpenAI; accessed 2026-10) was used interactively during subsequent mathematical exploration, assisting with the development, testing, and refinement of parts of the construction and proof. The author independently verified the final arguments and take full responsibility for the content of the paper.}

\vskip0.4cm
\noindent


%

\end{document}